\documentclass{article}
\usepackage{latexsym,amssymb,amsmath}
\usepackage{amsthm}
\usepackage{amsfonts}
\newtheorem{theorem}{Theorem}[section]

\newtheorem{note}[theorem]{Note}
\newtheorem{proposition}[theorem]{Proposition}
\date{}
\title
 {On a class of determinants}
\author{Milan Janji\'c \footnote{Department of Mathematics and Informatics, University of Banja Luka, Republic of Srpska}}

\begin{document}

\maketitle

\begin{abstract}
A class of determinants is introduced.
Different kind of mathematical objects, such as
Fibonacci, Lucas,  Tchebychev, Hermite, Laguerre, Legendre
polynomials, sums and covergents are represented as determinants from this class.

A closed formula in which arbitrary term of a homogenous linear
recurrence equation  is expressed in terms of the initial
conditions and the coefficients is proved.

\end{abstract}

\section{Introduction}
 As a particular case of Theorem 2.1 which is proved  in \cite{ja} we have the following theorem.
\begin{theorem}\label{tn1}
Let  $a_1,a_2,\ldots$ be a sequence which terms are from a commutative ring  $R.$ Assume that
\begin{equation}\label{fn}a_{k+1}=\sum_{i=1}^kp_{k,i}a_i,\;(p_{k,i}\in R, i=1,\ldots,k,\;k=1,\ldots.)\end{equation}
Then
\[a_{k+1}=a_1\left|\begin{array}{llllll}
p_{1,1}&p_{2,1}&p_{3,1}&\cdots&p_{k-1,1}&p_{k,1}\\
-1&p_{2,2}&p_{3,2}&\cdots&p_{k-1,2}&p_{k,2}\\
0&-1&p_{3,3}&\cdots&p_{k-1,3}&p_{k,3}\\
\vdots&\vdots&\vdots&\ddots&\vdots&\vdots\\
0&0&0&\cdots&p_{k-1,k-1}&p_{k-1,k}\\
0&0&0&\cdots&-1&p_{k,k}\end{array}\right|.\]
\end{theorem}

\section{Some consequences}
Several important mathematical objects may be represented in the form of a determinant of the above form.

 First of all, it is the case with the sequence of natural numbers, which may be defined by
the following recurrence relation:
\[a_1=1,\; a_{1+k}=a_1+a_k,\;(k=1,2,\ldots).\]
We thus obtain the following proposition.
\begin{proposition} Let $n$ be arbitrary natural number. Then
\[n=\left|\begin{array}{rrrrrr}
1&1&1&\cdots&1&1\\
-1&1&0&\cdots&0&0\\
0&-1&1&\cdots&0&0\\
\vdots&\vdots&\vdots&\ddots&\vdots&\vdots\\
0&0&\vdots&\cdots&1&0\\
0&0&\vdots&\cdots&-1&1\end{array}\right|,\]
where the size of the determinant is $n.$
\end{proposition}
Polynomials are also formed by the rule (\ref{fn}).
Namely, take a sequence $p_0,p_1,\ldots$ of elements of $R.$ Define coefficients in (\ref{fn}) in the following way.
\[p_{k,1}=p_{k-1},\;p_{k,k-1}=x,\;p_{k,i}=0.\;(i\not=1,k-1).\]
Then the formula (\ref{fn}) becomes
\[1,p_0,\;a_{1+k}=p_{k-1}a_1+xa_{k-1},\;(k=2,\ldots,n)\]
and obviously $a_{n+1}=f_n(x)=p_0x^n+p_1x^{n-1}+\cdots+p_0.$
We thus obtain the following propsition.
\begin{proposition} Let $f_n(x)=p_0x^n+p_1x^{n-1}+\cdots+p_0$ be a polynomial from $R[x].$ Then
\[f_n(x)=\left|\begin{array}{rrrrrr}
p_0&p_{1}&p_{2}&\cdots&p_{n-1}&p_n\\
-1&x&0&\cdots&0&0\\
0&-1&x&\cdots&0&0\\
\vdots&\vdots&\vdots&\ddots&\vdots&\vdots\\
0&0&\vdots&\cdots&x&0\\
0&0&\vdots&\cdots&-1&x\end{array}\right|.\]
\end{proposition}

Partial sums of a series may also be represented as determinants.
Taking, in particular, $x=1$ in the preceding equation we obtain $f_n(1)=\sum_{i=0}^np_i,$  that is,
\[\sum_{i=0}^np_i=\left|\begin{array}{rrrrrr}
p_0&p_{1}&p_{2}&\cdots&p_{n-1}&p_n\\
-1&1&0&\cdots&0&0\\
0&-1&1&\cdots&0&0\\
\vdots&\vdots&\vdots&\ddots&\vdots&\vdots\\
0&0&\vdots&\cdots&1&0\\
0&0&\vdots&\cdots&-1&1\end{array}\right|.\]

Some important classes of polynomials are given by recurrence relations. Applying Theorem \ref{tn1} we  may represent them as tridiagonal determinants.
This representation seems to be different from well-known  representation of orthonormal polynomials by Jacobi determinants.

Taking \[a_1=1,a_2=x,\;a_{k+1}=a_{k-1}+xa_k,\] Theorem \ref{tn1} yields the following proposition.
\begin{proposition}
If $F_1(x)=1,\;F_2(x)=x,\;F_3(x),\ldots$ are Fibonacci polynomials then
\begin{equation}\label{fibp}
F_{n+1}(x)=
\left|\begin{array}{rrrrrr}
x&1&0&\cdots&0&0\\
-1&x&1&\cdots&0&0\\
0&-1&x&\cdots&0&0\\
\vdots&\vdots&\vdots&\ddots&\vdots&\vdots\\
0&0&0&\cdots&x&1\\
0&0&0&\cdots&-1&x\end{array}\right|.\end{equation} The size of the determinant  is $n.$
\end{proposition}
Taking additionally $x=1$ we obtain the well-known formula for Fibonacci numbers.

\begin{equation}\label{fibb}
F_{n+1}=
\left|\begin{array}{rrrrrr}
1&1&0&\cdots&0&0\\
-1&1&1&\cdots&0&0\\
0&-1&1&\cdots&0&0\\
\vdots&\vdots&\vdots&\ddots&\vdots&\vdots\\
0&0&0&\cdots&1&1\\
0&0&0&\cdots&-1&1\end{array}\right|,\end{equation}
with $F_1=F_2=1.$

Under the conditions
\[a_1=1,a_2=x,\;a_3=2a_1+xa_2,\;a_{k+1}=a_{k-1}+xa_k,\;(k>2),\]
Theorem \ref{tn1} produces the following equation for  Lucas polynomials $L_n(x).$
\[L_{n+1}(x)=
\left|\begin{array}{rrrrrr}
x&2&0&\cdots&0&0\\
-1&x&1&\cdots&0&0\\
0&-1&x&\cdots&0&0\\
\vdots&\vdots&\vdots&\ddots&\vdots&\vdots\\
0&0&0&\cdots&x&1\\
0&0&0&\cdots&-1&x\end{array}\right|.\] The size of the determinant  is $n.$

The well-known recurrence relation for  Tchebychev polynomials  $T_k(x)$ of the first kind is:
\[T_0(x)=1,\;T_1(x)=x,\; T_{k}(x)=-T_{k-2}(x)+2xT_{k-1}(x),\;(k>2).\]
Theorem \ref{tn1} implies the following proposition.
\begin{proposition} For Tchebychev polynomials of the first kind $T_k(x)$  we have
\[T_{k}(x)=\left|\begin{array}{rrrrrr}
x&-1&0&\cdots&0&0\\
-1&2x&-1&\cdots&0&0\\
0&-1&2x&\cdots&0&0\\
\vdots&\vdots&\vdots&\ddots&\vdots&\vdots\\
0&0&0&\cdots&2x&-1\\
0&0&0&\cdots&-1&2x\end{array}\right|.\]
\end{proposition}
In the same way we obtain the following proposition.
\begin{proposition} For Tchebychev polynomials of the second kind $U_k(x)$  we have
\[U_{k}(x)=\left|\begin{array}{rrrrrr}
2x&-1&0&\cdots&0&0\\
-1&2x&-1&\cdots&0&0\\
0&-1&2x&\cdots&0&0\\
\vdots&\vdots&\vdots&\ddots&\vdots&\vdots\\
0&0&0&\cdots&2x&-1\\
0&0&0&\cdots&-1&2x\end{array}\right|.\]
\end{proposition}
\begin{note}
The preceding classes of polynomials  satisfy linear homogenous recurrence equation with constant coefficients.
\end{note}
 The next classes also satisfy  homogenous linear recurrence equations, but not with constant coefficients.

For Hermite polynomials $H_n(x)$ we have the  following recurrence relation.
\[H_0(x)=1,\;H_1(x)=2x,\; H_{n+1}(x)=-2nH_{n-1}(x)+2xH_{n}(x),\;(n\geq 2).\]
Applying Theorem \ref{tn1} for $n=1,2,\ldots$ we obtain the following proposition.
\begin{proposition} For Hermite polynomials  $H_n(x)$  we have
\begin{equation}\label{herp}
H_{n}(x)=
\left|\begin{array}{cccccc}
2x&-2&0&\cdots&0&0\\
-1&2x&-4&\cdots&0&0\\
0&-1&2x&\cdots&0&0\\
\vdots&\vdots&\vdots&\ddots&\vdots&\vdots\\
0&0&0&\cdots&2x&-2(n-1)\\
0&0&0&\cdots&-1&2x\end{array}\right|.\end{equation}
\end{proposition}
The recurrence relation for Legendre polynomials is:
\[P_0(x)=1,\;P_1(x)=x,\; P_{n+1}(x)=-\frac{n}{n+1}P_{n-1}(x)+\frac{2n+1}{n+1}xP_{n}(x),\;(n\geq 2).\]
Hence, the following proposition holds.
\begin{proposition} For Legendre polynomials  $P_k(x)$  we have
\[P_{n}(x)=
\left|\begin{array}{cccccc}
x&-\frac 12\vspace{2mm}&0&\cdots&0&0\\
-1&\frac 32 x&-\frac 23\vspace{2mm}&\cdots&0&0\\
0&-1&\frac53 x&\cdots&0&0\\
\vdots&\vdots&\vdots&\ddots&\vdots&\vdots\\
0&0&0&\cdots&\frac{2n-3}{n-1} x&-\frac{n-1}{n}\vspace{2mm}\\
0&0&0&\cdots&-1&\frac{2n-1}{n} x\vspace{2mm}\end{array}\right|.\]
\end{proposition}
In the same way we have the following proposition
\begin{proposition} For Laguerre polynomials we have
\[L_{n}(x)=\left|\begin{array}{cccccc}
1-x&-\frac 12\vspace{2mm}&0&\cdots&0&0\\
-1&\frac{3-x}{2}  &-\frac 23\vspace{2mm}&\cdots&0&0\\
0&-1&\frac{5-x}{3} x&\cdots&0&0\\
\vdots&\vdots&\vdots&\ddots&\vdots&\vdots\\
0&0&0&\cdots&-\frac{2n-3-x}{n-1}&-\frac{n-1}{n}\vspace{2mm}\\
0&0&0&\cdots&-1&\frac{2n-1-x}{n}\vspace{2mm}\end{array}\right|.\]
\end{proposition}

Taking  $p_{k,k}=p_k,\;p_{k,k-1}=1$ in (\ref{fn}) we obtain
\[a_1=1,p_1,\;a_{1+k}=a_{k-1}+p_ka_k,\;(k=2,\ldots).\]

Terms of this sequence are convergents of  $p_1,p_2,\ldots,p_n,\ldots$ and are denoted by $(p_1,p_2,\ldots,p_n).$
We thus have the following well-known formula.

\begin{equation}\label{kon}(p_1,p_2,\ldots,p_n)=\left|\begin{array}{llllll}
p_1&1&0&\cdots&0&0\\
-1&p_2&1&\cdots&0&0\\
\vdots&\vdots&\vdots&\vdots&\vdots&\vdots\\
0&0&0&\cdots&p_{n-1}&1\\
0&0&0&\cdots&-1&p_n\end{array}\right|.\end{equation}

Our next  result is a closed formula in which arbitrary term of a homogenous linear recurrence equation  is expressed in terms of initial conditions and coefficients of the equation.
Such one formula is usually obtained by the calculation of powers of a matrix, as in {\cite[p.\,62]{Age}}.

If $m$ is a fixed natural number, and $a(1),a(2),\ldots, a(m)$ arbitrary elements of $R.$
Consider a  homogenous linear recurrence equation of the form:
\begin{equation}\label{hlj}a(k)=p_1(k)a(k-m)+p_2(k)a(k-m-1)+\cdots+p_m(k)a(k-1),\;(k>m).\end{equation}

 According to Theorem \ref{tn1} we have the following theorem.
\begin{theorem}\label{linrek}
Let  $m$ be a fixed natural number,
  $a(1),a(2),\ldots,a(m)$ arbitrary elements of $R,$ and let $a(k),\;(k>m)$ be as in  (\ref{hlj}). Then
\[a(k)=\left|\begin{array}{ccccccccccc}
a(1)&a(2)&\cdots&a(m)&0&0&\cdots&0&\cdots&0&0\\
-1&0&\cdots&0&p_1(k)&0&\cdots&0&\cdots&0&0\\
0&-1&\cdots&0&p_2(k)&p_1(k)&\cdots&0&\cdots&0&0\\
\vdots&\vdots&\ddots&\vdots&\vdots&\vdots&\cdots&\vdots&\cdots&\vdots&\vdots\\
0&0&\cdots&-1&p_m(k)&p_{m-1}(k)&\cdots&0&\cdots&0&0\\
0&0&\cdots&0&-1&p_m(k)&\cdots&0&\cdots&0&0\\
\vdots&\vdots&\ddots&\vdots&\vdots&\vdots&\ddots&\vdots&\cdots&\vdots&\vdots\\
0&0&\cdots&0&0&0&\cdots&p_m(k)&\cdots&p_1(k)&0\\
0&0&\cdots&0&0&0&\cdots&-1&\cdots&p_2(k)&p_1(k)\\
\vdots&\vdots&\cdots&\vdots&\vdots&\vdots&\cdots&\vdots&\cdots&\vdots&\vdots\\
0&0&\cdots&0&0&0&\cdots&0&\cdots&p_m(k)&p_{m-1}(k)\\
0&0&\cdots&0&0&0&\cdots&0&\cdots&-1&p_m(k)
\end{array}\right|.\]
The size of the determinant is $k.$
\end{theorem}

We shall illustrate this formula by an example arising in differential equations.

Consider the equation
\[(x+1)y''+y'+xy=0,\;(y(0)=1,\;y'(0)=0).\]
Let $y(x)=\sum_{k=0}u(k)x^k$ be the solution.
 For coefficients $u(k)$ we easily obtain the following recurrence relation.
 \[u(k+2)=-\frac{k+1}{k+2}u(k+1)-\frac{1}{(k+1)(k+2)}u(k-1),\;(k\geq 0),\]
 with
 \[u(0)=1,\;u(1)=0,\; u(2)=0.\]

Applying Theorem \ref{linrek} for $u(k),\;(k>2)$ we obtain
\[u(k)=\left|\begin{array}{rrcccccrc}
1&0&1&0&0&\cdots&0\\
-1&0&0&-\frac{k-1}{k}\vspace{2mm}&0&\cdots&0\\
0&-1&0&0&-\frac{k-1}{k}\vspace{2mm}&\cdots&0\\
0&0&-1&-\frac{1}{(k-1)k}\vspace{2mm}&0&\cdots&0\\
0&0&0&-1&-\frac{1}{(k-1)k}\vspace{2mm}&\cdots&0\\
\vdots&\vdots&\vdots&\vdots&\vdots&\vdots&\vdots\\
0&0&0&0&0&0&-\frac{k-1}{k}\vspace{2mm}\\
0&0&0&0&0&0&0\\
0&0&0&0&0&-1&-\frac{1}{(k-1)k}\vspace{2mm}
\end{array}\right|.\]

\end{document}